\input epsf.tex 

\font\titlefont=cmr12 at 14pt
\font\sectionfont=cmbx10 at 12 pt
\font\curly=eusm10 at 10.7pt
\font\smallcurly=eusm10 at 9 pt
\font\Bbb=msbm10

\font\subsectfont=cmbx10 at 11pt

\def\subsect#1\par{\bigskip \leftline{\subsectfont #1 } \nobreak\smallskip}

\def\Fact#1. #2\par{\medskip\leftskip=10pt\noindent
{\bf Fact #1.}\enskip #2 \par\medskip\leftskip=0pt\noindent}

\magnification=\magstep1
\baselineskip=15pt
\mathsurround=1pt
\abovedisplayskip=8pt plus 3pt minus5pt
\belowdisplayskip=8pt plus 3pt minus5pt
\fontdimen16\textfont2=2.5pt  
\fontdimen17\textfont2=2.5pt  
\fontdimen14\textfont2=4.5pt  
\fontdimen13\textfont2=4.5pt

\def\page{\vfill\break}
\def\pf #1. {\noindent {\bf #1}.\enskip}

\def\int{\hbox{\rm int}}
\def\qed{\quad\hfill \rlap{$\sqcup$}$\sqcap$\par\medskip\smallskip} 

\def\:{\,\colon} 

\def\scrptS{\lower.2pt\hbox{\curly S}}
\def\smallscrptS{\lower.6pt\hbox{\smallcurly S}}

\def\S{\scrptS}
\def\s{{\smallscrptS}}

\def\R{\hbox{\Bbb R}}
\def\Z{\hbox{\Bbb Z}}

\def\phi{\varphi}

\def\bdy{\partial}

\def\longdownarrow{\hbox{\lower5pt\hbox{\hskip1pt$\downarrow$}\hskip-3.9pt\raise3pt\hbox{$|$}}}

\def\arrowunder #1 {\buildrel #1 \over{\longrightarrow}}
\def\leftarrowunder #1 {\buildrel #1 \over{\longleftarrow}}
\def\longarrowunder #1 {\buildrel #1 \over{-\hskip-5pt\longrightarrow}}
\def\longleftarrowunder #1 {\buildrel #1 \over{\longleftarrow\hskip-5pt-}}

\outer\def\section#1\par{\vskip0pt plus.05\vsize\penalty-250\vskip0pt plus-.05\vsize 
\vskip.5cm
\leftline{\sectionfont#1}\medskip}

\vbox{}\vskip .1truein
\centerline{\titlefont The Kirby Torus Trick for Surfaces}

\vskip.1truein
\centerline{Allen Hatcher}
\vskip .5truein

\leftskip30pt \rightskip30pt\noindent
{\bf Abstract}.  The classical theorem that topological surfaces can be triangulated is proved using the torus trick of Kirby plus a few basic facts about smooth or PL surfaces.

\leftskip0pt \rightskip0pt

\vskip 16pt

In the late 1960s the breakthrough that led to the Kirby--Siebenmann classification of piecewise linear (PL) structures on high-dimensional topological manifolds was the Kirby torus trick, appearing first in his 1969 paper [K] with a later exposition in Section I.3 of the Kirby--Siebenmann book [KS].  It still seems quite amazing that a construction this simple could so easily convert difficult topological problems into much more manageable ones in the PL category which could then be solved using existing techniques.  

In this paper we apply a scaled-down version of the torus trick to prove the existence and uniqueness of smooth or PL structures on topological surfaces, and in particular that surfaces can be triangulated.  The fact that the torus trick in high dimensions can also be applied for surfaces is mentioned in a footnote in [KS] (page 15).  However, this low-dimensional case does not seem to have appeared in the literature, so the purpose of the present paper is to provide a source for this approach while also serving as an introduction to the torus trick in a simple, easily visualized setting.

The special feature of the proof of triangulability of surfaces using the torus trick is that almost no point set topology is needed.  This is replaced instead by a few basic facts about smooth or PL surfaces whose proofs use only smooth or PL techniques.  In particular the topological Sch\"onflies theorem that every simple closed curve in the plane is the boundary of an embedded disk, which is an ingredient for most other triangulability proofs such as in [DM] or [T], is not needed here.  An exception is the proof by Moise in [M] where the triangulability of surfaces is proved before the Sch\"onflies theorem.  The first proof that surfaces can be triangulated is generally attributed to Rad\'o in 1925 [R].  That proof used the Sch\"onflies theorem without explicitly saying so, but the later exposition of Rad\'o's proof in [AS] makes the dependence clear.  

\section 1. The Theorems

Throughout the paper we take the term ``surface" to mean a $2$-dimensional manifold, possibly with boundary.  We assume manifolds are Hausdorff and are covered by a countable number of coordinate charts.  In dimension two it is well known and not hard to prove that the smooth and PL categories are equivalent.  For convenience we will work with smoothings rather than PL structures, but the reader who prefers PL structures could use these instead.

Here are the two theorems that will be proved using the torus trick:

\proclaim Theorem A. Every topological surface has a smooth structure. 

\proclaim Theorem B.  Every homeomorphism between smooth surfaces is isotopic to a diffeomorphism. 

Recall that a smooth structure on a surface $S$ is defined by a cover of $S$ by open sets $V_i$ together with homeomorphisms $h_i \: \R^2 \to V_i$ such that the compositions $h_i^{-1} h_j$ are diffeomorphisms $h_j^{-1}(V_i\cap V_j) \to h_i^{-1}(V_i\cap V_j)$ between open sets in $\R^2$. Thus the homeomorphisms $h_i$ provide local coordinate charts for $S$ which overlap in a smooth fashion.  Strictly speaking, a smooth structure is an equivalence class of collections $\{V_i,h_i\}$ as above, where two such collections are equivalent if their union is such a collection.

Theorem B implies that smooth structures are unique not just up to diffeomorphism, but up to isotopy.  This follows by applying the theorem to the identity map between two different smooth structures on the same surface.

The proof of Theorem A will show that a smooth structure exists extending a given smooth structure on a specified open set.  Similarly, Theorem B can be refined to the statement that if a homeomorphism is already a diffeomorphism near the boundary of a surface, the isotopy can be chosen to be fixed in a neighborhood of the boundary.  The proof will also show that the isotopy can be chosen to be arbitrarily small, so homeomorphisms can be approximated arbitrarily closely by diffeomorphisms.

The fact that homeomorphisms of surfaces are isotopic to PL homeomorphisms and hence to diffeomorphisms was proved by Epstein in 1966 [E], although at least for closed orientable surfaces this can be derived from results of Baer in 1928 [B] on mapping class groups.

\medskip

Let us say something about the basic ideas underlying the proofs of Theorems A and B using the torus trick.  For both theorems one would like to isotope homeomorphisms to diffeomorphisms.  This is evident for Theorem~B, while for 
Theorem~A one would like to isotope a system of coordinate charts $h_i$ so that the coordinate-change homeomorphisms $h_i^{-1} h_j$ of open sets in $\R^2$ are diffeomorphisms.  In both cases one is working with homeomorphisms of smooth surfaces.  Smooth surfaces have smooth triangulations, so one can proceed inductively, first smoothing a homeomorphism near vertices, then extending this to a smoothing near edges, then extending over the 2-simplices.  It turns out that the torus trick is needed only for smoothing near vertices, the other two cases being more straightforward.  
A neighborhood of a vertex is a copy of $\R^2$ and a homeomorphism takes this onto another copy of $\R^2$. Thus one has a homeomorphism of $\R^2$ which one would like to deform by isotopy to a diffeomorphism in a neighborhood of a point, staying fixed outside a larger neighborhood of the point. 

A special type of homeomorphism of $\R^2$ which is more tractable than a general homeomorphism is one which is bounded, meaning that the distances that points are moved by the homeomorphism are bounded by a constant.  In this case if we regard $\R^2$ as the interior of the closed unit disk $D^2$ by radial reparametrization, then a bounded homeomorphism of $\R^2$ becomes a homeomorphism of the open disk that extends to a homeomorphism of $D^2$ which is the identity on $\bdy D^2$. Such a homeomorphism of $D^2$ is isotopic to the identity by the classical Alexander trick, the isotopy which at time $t$ consists of a shrunken copy of $h$ in the disk of radius $1-t$ centered at the origin, extended by the identity outside this disk.  Restricting to the open unit disk and regarding the open disk again as $\R^2$, this gives an isotopy of a bounded homeomorphism of $\R^2$ to the identity.

The Kirby torus trick is a method for reducing to the case of bounded homeomorphisms of $\R^2$.  This is motivated by the observation that the lift of a homeomorphism of the torus $T$ to the universal cover $\R^2$ of $T$ is a bounded homeomorphism, provided that the homeomorphism of $T$ induces the identity map on $\pi_1(T)$, which is easy to achieve by composing with a homeomorphism of $T$ defined by an element of $GL_2(\Z)$.  The goal of the torus trick is to convert a homeomorphism of $\R^2$ into a homeomorphism of $T$ which can then be lifted to a bounded homeomorphism of $\R^2$.  

This is done somewhat indirectly.  In outline:  A homeomorphism $h\:\R^2\to \R^2$ pulls the standard smooth structure on $\R^2$ back to another smooth structure on $\R^2$.  Since $T$ with a point deleted immerses in $\R^2$, the new smooth structure on $\R^2$ then induces a smooth structure on the punctured torus.  Basic facts about smooth surfaces imply first that this new smooth structure on the punctured torus can be extended to a smooth structure on $T$, and then that this new smooth structure on $T$ is diffeomorphic to the standard smooth structure.  Thus we have a diffeomorphism from $T$ with its new smooth structure to $T$ with its standard smooth structure. After composing this diffeomorphism of $T$ with a linear diffeomorphism to make it induce the identity on $\pi_1$ we obtain a diffeomorphism between the two smooth structures on $T$ that lifts to a bounded homeomorphism of $\R^2$.  This gives a homeomorphism of an open disk which extends by the identity to a homeomorphism of $\R^2$. The Alexander trick gives an isotopy $G_t\:\R^2\to\R^2$ from the identity to this  homeomorphism of $\R^2$. One then checks that the desired isotopy of the original homeomorphism $h$ is given by the composition $hG^{-1}_t$.

\smallskip

The torus trick was originally applied for $n$-dimensional manifolds with $n\geq 5$.  The input needed was information about the different smooth or PL structures on the $n$-torus $T^n$ and more generally products $D^k \times T^{n-k}$.  The difference between smooth and PL structures was fairly well understood when the torus trick was discovered, with exotic spheres providing obstructions to the existence and uniqueness of smooth structures on PL manifolds, so the remaining question was the existence and uniqueness of PL structures on topological manifolds.  It turned out that only in the case of $D^3\times T^{n-3}$ was there a nontrivial obstruction, and this was only a mod 2 obstruction.  From this Kirby and Siebenmann concluded that there was a single well-defined obstruction in $H^4(M;\Z_2)$ to the existence of a PL structure on a topological $n$-manifold $M$ when $n\geq 5$, with isotopy classes of PL structures on $M$ corresponding bijectively to elements of $H^3(M;\Z_2)$ when a PL structure exists.  

For manifolds of dimension four much less is known about the existence and uniqueness of PL structures. In dimensions less than seven PL and smooth structures are equivalent, so the question is usually studied as existence and uniqueness of smooth structures on $4$-manifolds.  From the current limited state of knowledge it is apparent that there is a much wider chasm between PL and topological $4$-manifolds than in higher dimensions, and there is no hope of the torus trick leading to an obstruction theory for existence and uniqueness of PL structures analogous to the theory in high dimensions.

For 3-manifolds the torus trick can be made to work.  An extra complication not present in the surface case is avoiding the smooth Poincar\'e conjecture.  There is a 1974 paper by A.\hskip1pt J.\hskip1pt S.\hskip1.5pt  Hamilton [H1] that presents this approach to triangulating $3$-manifolds.  An alternative version of the argument will be given in [H3], filling what appears to be a gap in Hamilton's proof.

\section 2. Reduction to Handle Smoothing.

Theorems A and B will be deduced from a result about smoothing handles in surfaces.  An $i$-handle in an $n$-manifold is usually taken to be a product $D^i \times D^{n-i}$ of closed disks but it will be more convenient here to use a product $D^i \times \R^{n-i}$ instead.  For handles in surfaces there are three cases: $0$-handles, $1$-handles, and $2$-handles.  The proofs that handles can be smoothed are somewhat different in each case, with the torus trick being needed only for $0$-handles.  The three cases could be combined into a single statement, but here are the separate statements:

\proclaim Handle Smoothing Theorem. Let $S$ be a smooth surface. Then{\rm\hskip1pt :}\hfill\break
\hbox{\rm ---}\enskip A topological embedding $\R^2\to S$ can be isotoped to be a smooth embedding in a neighborhood of the origin, staying fixed outside a larger neighborhood of the origin. \hfill\break
\hbox{\rm ---}\enskip A topological embedding $D^1\times \R\to S$ which is a smooth embedding near $\bdy D^1\times \R$ can be isotoped to be a smooth embedding in a neighborhood of $D^1\times 0$, staying fixed outside a larger neighborhood of $D^1\times 0$ and in a neighborhood of $\bdy D^1\times \R$. \hfill\break
\hbox{\rm ---}\enskip A topological embedding $D^2\to S$ which is a smooth embedding in a neighborhood of $\bdy D^2$ can be isotoped to be a smooth embedding on all of $D^2$, staying fixed in a smaller neighborhood of $\bdy D^2$.

\smallskip
It should be remembered that a smooth embedding is more than just a topological embedding which is a smooth map since the differential must also be nonsingular at each point.

Assuming the Handle Smoothing Theorem, let us now deduce Theorems A and B.  The proofs will use a few standard results about smooth surfaces, and these will be stated just as facts, with proofs given (or sketched) in the last section of the paper for those who would like to see them.  The same procedure will be followed later when we prove the Handle Smoothing Theorem.

\medskip
\pf Proof of Theorem A. For a surface $S$ without boundary, choose a system of coordinate charts $h_i \: \R^2\to S$, $i=1,2,\cdots$.  Let $V_i$ be the image $h_i(\R^2)$, an open set in $S$.  Each $V_i$ has a smooth structure induced by $h_i$ from the standard smooth structure on $\R^2$.  The problem is to  make these smooth structures agree on the overlaps of different $V_i$'s. 

Let $U_n$ be the union $V_1 \cup \cdots \cup V_n$ and assume by induction on $n$ that we have already modified $h_1,\cdots,h_n$ by isotopy to give a well-defined smooth structure on $U_n$.  The induction starts with the case $n=1$ where no modification of $h_1$ is needed for it to give a smooth structure on $U_1=V_1$.  The induction step will be to isotope $h_{n+1}$ so that it restricts to a diffeomorphism from the open set $W=h^{-1}_{n+1}(U_n\cap V_{n+1})$ in $\R^2$ onto $U_n\cap V_{n+1}$ where $U_n\cap V_{n+1}$ has the smooth structure obtained by restriction of the smooth structure constructed inductively on $U_n$.  After this isotopy of $h_{n+1}$ we will then have a well-defined smooth structure on $U_{n+1}$.  These inductively constructed smooth structures on all the sets $U_n$ then give a smooth structure on $S$.

\page

To construct the isotopy of $h_{n+1}$ we will use the following general statement:  

\Fact 1. An open set $W\subset\R^2$ can be triangulated so that the size of the simplices approaches $0$ at the frontier of $W$.  

Having such a triangulation, the idea will be to isotope $h_{n+1}$ to make it a diffeomorphism on $W=h^{-1}_{n+1}(U_n\cap V_{n+1})$ inductively over neighborhoods of the skeleta of the triangulation.  First, we choose disjoint $0$-handles about all the vertices, with the vertices located at the centers of the $0$-handles.  Applying $0$-handle smoothing in each of these handles then gives an isotopy of the restriction of $h_{n+1}$ to $W$ to a new $h_{n+1}$ which is a smooth embedding of an open neighborhood $N_0$ of the $0$-skeleton into $U_n$.  Since the simplices become small near the frontier of $W$, so do the $0$-handles, so this isotopy of $h_{n+1}$ on $W$  can be extended by the constant isotopy on $\R^2 - W$.  

Next, for each edge of the triangulation choose a $1$-handle $D^1\times \R$ with ${D^1\times 0}$ contained in the edge and $\bdy D^1 \times \R$ contained in $N_0$.  We can assume all these $1$-handles are disjoint. Applying $1$-handle smoothing to all the $1$-handles gives an isotopy of $h_{n+1}$ as before so that it becomes a smooth embedding of an open neighborhood $N_1$ of the $1$-skeleton.  Finally, we apply $2$-handle smoothing to isotope $h_{n+1}$ further to be a smooth embedding on the $2$-handles obtained by slightly shrinking each $2$-simplex of the triangulation, with the boundary of the $2$-handle in $N_1$. The resulting $h_{n+1}$ is a homeomorphism whose restriction  $W\to U_n\cap V_{n+1}$ is locally a diffeomorphism and hence globally a diffeomorphism.  This finishes the induction step.  After all the coordinate mappings $h_n$ have been isotoped we then have a smooth structure on all of $S$, finishing the proof of Theorem~A when $\bdy S$ is empty. 

To treat the case that $\bdy S$ is nonempty we will use the fact that $\bdy S$ has a collar neighborhood in $S$.  This is a general fact about topological manifolds, with an elementary proof due to Connelly [C] that uses a partition of unity to piece together local collars. (For compact manifolds this argument is given in Proposition~3.42 in [H2] and the argument extends easily to the noncompact case.)  An open collar neighborhood $U_0$ of $\bdy S$ has a smooth structure since $1$-manifolds are smoothable.  Then the inductive procedure above extends this to a smooth structure on all of $S$ using coordinate charts $h_i \:  \R^2 \to S$ covering the interior of $S$, with $U_n=U_0 \cup V_1 \cup\cdots \cup V_n$. 
\qed

\medskip
\pf Proof of Theorem B. Let $f \:  S\to S'$ be a homeomorphism of smooth surfaces.  We will isotope $f$ to be a diffeomorphism by  isotoping it on one handle at a time, with handles obtained from a triangulation of $S$ using the following:  

\page

\Fact 2. A smooth surface $S$ has a smooth triangulation, with $\bdy S$ a subcomplex.

To apply this, suppose first that $\bdy S$ is empty.  As in the proof of Theorem~A we can do $0$-handle smoothing to isotope $f$ to be a diffeomorphism near all vertices of the triangulation of $S$, then apply $1$-handle smoothing to isotope the new $f$ to be a diffeomorphism near the $1$-skeleton of the triangulation, and finally apply $2$-handle smoothing to make the resulting $f$ a diffeomorphism on all of $S$. 

When $\bdy S$ is nonempty we first choose a smooth collar on $\bdy S$ and isotope $f$ to preserve the collar parameter in a smaller collar.  The restriction of $f$ to $\bdy S$ is isotopic to a diffeomorphism, and this isotopy can be extended to an isotopy of $f$ which is constant outside the collar and ends with a new $f$ which is a diffeomorphism in a smaller collar.  This gives a reduction to the case that $f$ is already smooth on a  neighborhood of $\bdy S$.  Then we can apply the smoothing procedure from the preceding paragraph to make $f$ a diffeomorphism on the handles corresponding to simplices not contained in $\bdy S$.
\qed

By choosing the triangulation of $S$ to have small simplices we can make the isotopy of $f$ to a diffeomorphism small as well. In particular, homeomorphisms can be approximated by diffeomorphisms.  

For the proofs of Theorems~A and~B it would suffice to prove handle smoothing in the special case that the target surface $S$ is $\R^2$ since for Theorem~A the proof only involved handles contained in coordinate charts, while for Theorem~B a triangulation of $S$ can be chosen with each simplex contained in a coordinate chart.  Restricting $S$ to be $\R^2$ does not seem to offer any significant simplifications in the proof of handle smoothing, however.

The deduction of Theorems~A and~B from the Handle Smoothing Theorem used arguments that apply not just for surfaces but for manifolds of any dimension in which all handles can be smoothed.  However, the only higher dimension where this happens is dimension three.

\section 3. The Torus Trick:  Smoothing 0-Handles.

\parshape=4 0pt\hsize 0pt\hsize 0pt.75\hsize  0pt.75\hsize 
This is the hardest case, where the torus trick is used.  We view the torus $T$ as the orbit space $\R^2/\Z^2$, with a basepoint $0$ that is the image of $0\in\R^2$.  Deleting some other point $*\in T$ yields a punctured torus $T'$.  \vadjust{\hfill\smash{\lower73pt\llap{\epsfbox{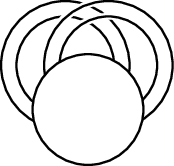}}}}\ignorespaces 
We can immerse $T'$ in $\R^2$, with $0\in T'$ mapping to $0\in\R^2$, by viewing $T'$ as the interior of the surface obtained from a disk by attaching two $1$-handles, then letting the immersion be an embedding of the disk, with the two $1$-handles individually embedded so that their images cross in $\R^2$ as shown in the figure. 

\page

Let $h\:\R^2\to S$ be a topological embedding into the smooth surface $S$.  Via $h$, the smooth structure on $S$ induces a smooth structure $\S$ on $\R^2$.  This is likely to be quite different from the usual smooth structure if $h$ is far from being a smooth embedding.  For example, smooth curves in $S$ that are in the image of $h$ pull back via $h^{-1}$ to curves in $\R^2$ that can be quite wild when looked at with the naked eye but are in fact smooth in the structure $\S$.  

The immersion $T'\to\R^2$ pulls back the smooth structure $\S$ on $\R^2$ to a smooth structure on $T'$ that we denote $T'_{\s}$.

\Fact 3.  For each smooth structure on the punctured torus $T'$ there is a compact set in $T'$ whose complement is diffeomorphic to $S^1\times\R$.

\noindent
The complement of the compact set is a punctured neighborhood of the puncture point in $T$, so Fact~3 allows the smooth structure $T'_{\s}$ to be extended to a smooth structure $T_{\s}$ on~$T$.   

\Fact 4.  Every smooth structure on a torus $S^1\times S^1$ is diffeomorphic to the standard smooth structure.

Thus there is a diffeomorphism $g\:T_{\s}\to T$.  Note that $g$, viewed just as a homeomorphism from $T$ to $T$, is likely to be as complicated as the original $h$ locally.

We can normalize $g$ so that it takes $0$ to $0$ by composing with rotations in the $S^1$ factors of the target if necessary.  We can then further normalize so that $g$ induces the identity automorphism of $\pi_1(T,0)=\Z^2$ by composing with a diffeomorphism in the target given by a suitable element of $GL_2(\Z)$ acting on $\R^2/\Z^2$. After this normalization, $g$ lifts to the universal covers as a diffeomorphism $\widetilde g\:\R^2_{\s}\to \R^2$ fixing $\Z^2$, where the subscript $\S$ denotes the smooth structure lifted from $T_{\s}$.  The lift $\widetilde g $ is periodic, satisfying $\widetilde g(x+m,y+n)=\widetilde g(x,y)+(m,n)$ for all $(x,y)\in \R^2$ and $(m,n)\in \Z^2$. 

The key point of these constructions is that $\widetilde g$ is a bounded homeomorphism $\R^2\to\R^2$, so the distance $|\,\widetilde g(x,y) - (x,y)|$ is less than a fixed constant for all $(x,y)\in\R^2$.  This is certainly true as $(x,y)$ varies over the square $I^2$ by compactness, and then it holds over all of $\R^2$ by periodicity.  

If we identify $\R^2$ with the interior of $D^2$ by a radial reparametrization, then $\widetilde g$ becomes a homeomorphism of the interior of $D^2$ that extends via the identity on $\bdy D^2$ to a homeomorphism $G\:D^2\to D^2$, as a result of the boundedness condition.  This can be seen by considering polar coordinates $(r,\theta)$, since as a disk in $\R^2$ of fixed radius and varying center moves out to infinity, the variation in the $\theta$-coordinates of points in the disk approaches zero, and after the radial reparametrization to move $\R^2$ inside $D^2$, the variation in the $r$-coordinates of points in the disk also approaches zero.  We can choose the identification of $\R^2$ with the interior of $D^2$ to be the identity near $0$, and then $G=\widetilde g = g$ near $0$.  Finally, we can extend $G\:D^2\to D^2$ to a homeomorphism $G\:\R^2\to \R^2$ that is the identity outside $D^2$.  

By the Alexander trick $G$ is isotopic to the identity.  This isotopy can be obtained by replacing the unit disk $D^2$ in the preceding paragraph by disks of progressively smaller radius centered at $0$, limiting to radius zero.  Reversing the time parameter of this isotopy, we obtain an isotopy $G_t\:\R^2\to \R^2$ with $G_0$ the identity and $G_1=G$.
We can collect all the maps described above into one diagram:
$$\eqalign{
& \hskip5pt \R^2 \longleftarrow T' \subset\hskip1pt T \longleftarrow \R^2 \hskip3pt\subset\hskip3pt D^2 \hskip3pt\subset\hskip3pt \R^2 \cr
& h \longdownarrow \hskip 47pt g\longdownarrow \hskip19pt \widetilde g \longdownarrow \hskip16pt G_t \hskip0pt\longdownarrow \hskip15pt G_t \hskip0pt\longdownarrow \cr
& \hskip5pt S \hskip51pt T \longleftarrow \R^2 \hskip3pt\subset\hskip3pt D^2 \hskip3pt\subset\hskip3pt \R^2 \cr 
}$$
\vskip-3pt
\noindent

We check now that the desired isotopy of $h$ is given by $h_t=hG_t^{\hskip1pt -1}$.  We have $h_0=h$ since $G_0$ is the identity.  Also $h_t$ is stationary outside $D^2$ since $G_t$ is the identity there.  Since $G_t$ fixes the origin for all $t$, we have $h_t(0)=h(0)$ for all $t$.  
Finally, to check that $h_1$ is a smooth embedding near $0$ with respect to the standard smooth structure on $\R^2$ we have $h_1=hG_1^{\hskip1pt -1}=hG^{\hskip.7pt -1}$ with $G^{\hskip.7pt -1}$ a diffeomorphism from the standard smooth structure to the smooth structure $\S$ near $0$,  and $h$ carries the smooth structure $\S$ to the smooth structure that was given on $S$.  \qed

\section 4. Smoothing 1-Handles and 2-Handles.

For smoothing a 1-handle we are given an embedding $h\: D^1\times \R\to S$ which is smooth near $\bdy D^1\times \R$.  The smooth structure on $S$ pulls back via $h$ to a smooth structure $\S$ on $D^1\times \R$ which agrees with the standard structure near $\bdy D^1\times\R$.  

\Fact 5. Every smooth structure on $D^1\times\R$ that is standard near $\bdy D^1\times\R$ is diffeomorphic to the standard structure via a diffeomorphism that is the identity near $\bdy D^1\times\R$.

Thus there is a diffeomorphism $g\:(D^1\times \R)_{\s}\to D^1\times \R$ which is the identity near $\bdy D^1\times\R$.  

Since $g^{-1}(D^1\times D^1)$ is compact, there is a number $a > 1$ such that $g^{-1}(D^1\times D^1) \subset D^1\times (-a,a)$.  Choose an embedding $\sigma\:D^1\times \R\to D^1\times\R$ which is the identity on $D^1\times [-a,a]$ and has image $\bigl(D^1\times [-b,b]\bigr) - \bigl(0\times \{b,-b\}\bigr)$ for some $b > a$, as indicated in the figure below.  

\vskip3pt
\centerline{\epsfbox{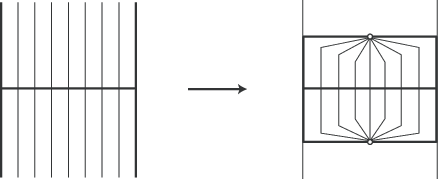}}

\noindent
The composition $\sigma g\hskip.5pt\sigma^{-1}$ is then a homeomorphism of the image of $\sigma$.  
This homeomorphism can be extended by the identity outside the image of $\sigma$ to yield a homeomorphism $G\:D^1\times\R\to D^1\times\R$.  If we perform the Alexander isotopy on $G$, using the rectangle $D^1 \times [-b,b]$ this time instead of the disk $D^2$, we obtain, after reversing the time parameter, an isotopy $G_t\: D^1\times\R\to D^1\times\R$ from the identity to $G$.  
This is stationary outside $D^1\times [-b,b]$, and it is also stationary near $\bdy D^1 \times \R$ since $g$ is the identity near $\bdy D^1 \times \R$. 

The composition $hG^{-1}_t$  is an isotopy $h_t \: D^1\times\R\to S$ from $h_0=h$ to $h_1= hG^{-1}$ that is stationary near $\bdy D^1 \times \R$ and outside $D^1\times [-b,b]$. On $D^1\times D^1$ we have $h_1=hg^{-1}$ since $G^{-1}=g^{-1}$ on $D^1\times D^1$.  The composition $hg^{-1}$ takes the standard smooth structure on $D^1\times \R$ first to the smooth structure $\S$ on $D^1\times \R$, then to the given smooth structure on~$S$.  Thus in a neighborhood of $D^1\times 0$ in $D^1\times \R$, $h_1$ is a smooth embedding with respect to the standard smooth structure on $D^1\times \R$.  The isotopy $h_t$ is therefore the desired isotopy of $h$, converting it into a smooth embedding near $D^1\times 0$, staying fixed near $\bdy D^1\times \R$ and outside $D^1\times [-b,b]$.  \qed

Finally we have the case of smoothing 2-handles.  Here we start with an embedding $h\: D^2\to S$ which is smooth near $\bdy D^2$.  The smooth structure on $S$ pulls back to a smooth structure $D^2_{\s}$ which is standard near $\bdy D^2$.  

\Fact 6. Every smooth structure on $D^2$ that is standard near $\bdy D^2$ is diffeomorphic to the standard structure via a diffeomorphism that is the identity near $\bdy D^2$.

Thus there is a diffeomorphism $g\:D^2_{\s}\to D^2$ which is the identity near $\bdy D^2$.  We can apply the Alexander trick directly to $g$, yielding an isotopy $g_t\: D^2\to D^2$ from the identity to $g$.  Then an isotopy from $h$ to a smooth embedding is given by $h_t=hg_t^{-1}$, staying fixed on a neighborhood of $\bdy D^2$.  \qed

\page

\section 5.  Proofs of the Facts.

\proclaim Fact 1.  An open set $W\subset\R^2$ can be triangulated so that the size of the simplices approaches $0$ at the frontier of $W$.  

\pf  Proof.  The horizontal and vertical lines in $\R^2$ through points in the integer lattice $\Z^2$ divide $\R^2$ into closed unit squares.  To start, choose all such squares that are contained in the open set $W$.  For the squares that are not contained in $W$, divide each of these into four squares with sides of length $1/2$ and take all the smaller squares that lie in $W$.  Now repeat this process indefinitely, at each stage subdividing the squares of the preceding stage not contained in $W$ into four subsquares and choosing the ones contained in $W$.  The union of all the squares chosen by this iterative process is $W$ since the distance from any point in $W$ to $\R^2-W$ is positive.  Only finitely many of the squares touch a given square since the distance from the square to $\R^2-W$ is positive. The chosen squares give a cellulation of $W$, and this can be subdivided to a triangulation by adding a new vertex at the center of each square and joining this vertex to the vertices in the boundary of the square. \qed

\proclaim Fact 2. Every smooth surface $S$ has a smooth triangulation with $\bdy S$ a subcomplex.

\pf Proof. 
Consider first the case that $\bdy S$ is empty.  We will construct a smooth cellulation of $S$ with $2$-cells that are polygons glued together along their edges.
Choose a Morse function $S\to \R$ that is proper (inverse images of compact sets are compact) and has all its critical points on distinct levels.  The existence of such a function is a standard fact about compact smooth manifolds, and the noncompact case can be deduced from the compact case in the following way.  First construct a proper smooth function $S\to [0,\infty)$, as for example in Proposition~2.28 in [L].  This function has a dense set of noncritical values by Sard's Theorem, and the preimages of a sequence of noncritical values approaching $\infty$ cut $S$ into compact subsurfaces.  Morse functions on these subsurfaces taking constant values on boundary components can then be combined to give a Morse function on $S$, and a small perturbation of this puts all critical points on distinct levels.  

Having a proper Morse function $S\to [0,\infty)$ with all its critical points on distinct levels, we can then cut $S$ along noncritical levels separating all critical levels to obtain a decomposition of $S$ into pieces that are diffeomorphic to disks, annuli, pairs of pants, or annuli with a 1-handle attached to one boundary circle so as to produce a nonorientable surface.  A surface of this last type can be further cut along a circle in its interior to produce a pair of pants and a M\"obius band.  Thus $S$ is decomposed into disks, annuli, pairs of pants, and M\"obius bands, glued together along their boundary circles.  From this a cellulation is easily obtained by inserting one vertex in each circle, then in each annulus connecting the two vertices in its boundary by an arc cutting the annulus into a square, in each pair of pants connecting the three vertices in its boundary by two arcs cutting the pair of pants into a polygon (a heptagon in fact), and in each M\"obius band inserting an edge to cut it into a triangle.  (Think of the M\"obius band as a rectangle with two opposite edges identified, then cut the rectangle along a diagonal before identifying the opposite edges.)  Having a smooth cellulation of $S$, we can then subdivide it to a smooth triangulation, although for the application of Fact~2 in Section~2 a smooth cellulation would suffice. 

The extension of this construction to the case that $\bdy S$ is nonempty is similar, using a proper Morse function with all its critical points in the interior of $S$ and with its restriction to $\bdy S$ a proper Morse function on $\bdy S$.  Details are left to the reader. 
\qed

\proclaim Fact 3. If $T'_{\s}$ is a smooth structure on the punctured torus $T'$ then there is a compact subsurface of $T'_\s$ whose complement is diffeomorphic to $S^1\times\R$.

\pf Proof.  As in the proof of Fact 2, $T'_{\s}$ can be cut along a collection of disjoint smooth circles $C_i$ to produce pieces $P_j$ that are disks, annuli, pairs of pants, and M\"obius bands, but the last of these cannot occur since $T'_{\s}$ is orientable.  The pattern in which the pieces $P_j$ are assembled to form $T'_{\s}$ is described by a graph $G$ having one vertex for each $P_j$ and one edge for each circle $C_i$.  There is a quotient map $q\:T'_{\s}\to G$ projecting a product neighborhood of each $C_i$ onto its interval factor, which becomes the corresponding edge of $G$, and collapsing the complementary components of these annuli to points, the vertices of $G$.  The induced homomorphism $q_*\:\pi_1T'_{\s}\to\pi_1G$ is surjective and in fact split since $q$ has a section up to homotopy.  Thus the free group $\pi_1G$ is a quotient of $\pi_1T'_{\s}$, so it is finitely generated.  This implies that there is a finite subgraph $G_0\subset G$ such that the closure of $G-G_0$ consists of a finite number of trees.  Only one of these trees can be noncompact since $T'_{\s}$ has only one end.  The one-endedness also implies that this noncompact tree consists of an infinite subtree homeomorphic to $[0,\infty)$ with finite subtrees attached to it.  We can eliminate these finite subtrees by deleting the edges leading to vertices corresponding to disk pieces $P_j$, discarding the circles $C_i$ corresponding to these edges.  After this has been done the end of $T'_{\s}$ consists of an infinite sequence of annuli glued together to form an infinite cylinder.  The glueings are smooth so the cylinder is diffeomorphic to a standard cylinder.  \qed

\page

\proclaim Fact 4.  Every smooth structure on a torus $S^1\times S^1$ is diffeomorphic to the standard smooth structure.

This of course follows from the classification of smooth closed surfaces, but let us give an argument similar the one for Fact~3.

\medskip

\pf Proof.  Let the given smooth structure on the torus be denoted  $T_{\s}$ and choose a Morse function $T_{\s} \to \R$. The associated graph $G$ is now a finite graph with $\pi_1G$ a quotient of the abelian group $\pi_1T_{\s}$, so $\pi_1G$ is either $\Z$ or the trivial group.  In the former case we can reduce $G$ to a circle by eliminating vertices corresponding to disk pieces as before, and then we see that $T_{\s}$ consists of annuli glued together in a cyclic pattern, so $T_{\s}$ is diffeomorphic to a torus or a Klein bottle, but the latter is ruled out since $\pi_1T_{\s}$ is abelian (or since $T_{\s}$ is orientable).  In the other case that $\pi_1 G$ is trivial we can reduce $G$ to a single edge, which would mean that $T_{\s}$ was a sphere, so this case cannot occur.  \qed

\proclaim Fact 5. Every smooth structure on $D^1\times\R$ that is standard near $\bdy D^1\times\R$ is diffeomorphic to the standard structure via a diffeomorphism that is the identity near $\bdy D^1\times\R$.

\pf Proof. Let $(D^1\times \R)_{\s}$ be the given smooth structure on $D^1\times \R$. The projection $D^1\times\R\to\R$ can be perturbed to a proper Morse function $f\:(D^1\times\R)_{\s}\to\R$ staying fixed near $\bdy D^1\times\R$ where it is already smooth without critical points.  We can assume all the critical points of $f$ lie in distinct levels.  Each noncritical level consists of a single arc and finitely many circles.  Cutting along these curves for a set of noncritical levels separating all the critical levels produces pieces that are disks, annuli, or pairs of pants as in the earlier cases, and there are now also two new kinds of pieces: rectangles, and rectangles with an open disk removed.  The associated graph $G$ is a tree since $\pi_1(D^1\times\R)=0$, and $G$ has two ends since $(D^1\times\R)_{\s}$ has two ends.  Thus $G$ consists of a subgraph homeomorphic to $\R$ with finite trees attached to it.  These finite trees can be eliminated as before, and $f$ can be modified to a proper Morse function without critical points, equal to the old $f$ near $\bdy D^1\times\R$.  
The new $f$ is then the second coordinate of a diffeomorphism $g\:(D^1\times\R)_{\s}\to D^1\times\R$ 
whose first coordinate is obtained by flowing up or down from $f^{-1}(0)$ along the integral curves of the gradient vector field of $f$.  These curves are standard near $(\bdy D^1\times\R)_{\s}$ so $g$ can be taken to be the identity there.   \qed

\proclaim Fact 6. Every smooth structure on $D^2$ that is standard near $\bdy D^2$ is diffeomorphic to the standard structure via a diffeomorphism that is the identity near $\bdy D^2$.

\pf Proof.  Here we extend the radial coordinate in $D^2$ near $\bdy D^2$ to a Morse function $f\:D^2_{\s}\to[0,1]$ with $f^{-1}(1)=\bdy D^2_{\s}$ and all critical points in the interior of $D^2_{\s}$.  The associated graph is a finite tree, and it can be simplified as before until it is a line segment, with a new $f$ having a single critical point, a local minimum.  Then $f$ together with the flow lines of its gradient field can be used to construct a diffeomorphism $g\:D^2_{\s}\to D^2$ which is the identity near $\bdy D^2_{\s}$. \qed

A remark on the preceding proof:  A diffeomorphism $g\:D^2_{\s}\to D^2$ that is not the identity near $\bdy D^2$ can be modified to one that is by first taking the restriction of $g$ to $\bdy D^2$ and extending this to a diffeomorphism $\widehat g$ of $D^2$, then replacing $g$ by $\widehat g\hskip1pt^{-1} g$ to get a diffeomorphism $D^2_{\s}\to D^2$ that is the identity on $\bdy D^2$.  This can then be isotoped to be the identity in a neighborhood of $\bdy D^2$.  Similarly in Fact~5 a diffeomorphism $( D^1\times\R)_{\s}\to D^1\times \R$ can be modified to be the identity near $\bdy D^1\times\R$.

\medskip

\leftskip 20pt
\parindent=-20pt

\vskip-15pt

\section References

[AS] L.\hskip1pt V.\hskip2pt Ahlfors and L.\hskip2pt Sario, {\it Riemann Surfaces}, 1960.

[B] R.\hskip2pt Baer, Isotopie von Kurven auf orientierbaren, geschlossenen Fl\"achen und ihr Zusammenhang mit der topologischen Deformation der Fl\"achen, {\it J.\ Reine u.\ Angew.\ Math.\ }159 (1928), 101--116.

[C] R.\hskip2pt Connelly, A new proof of Brown's collaring theorem, {\it Proc.\ A.M.S.\ }27 (1971), 180--182.

[DM] P.\hskip1pt H.\hskip2pt Doyle and D.\hskip1pt A.\hskip2pt Moran,  A short proof that compact 2-manifolds can be triangulated, {\it Inventiones math.\ }5 (1968), 160\hskip1pt--162.

[E] D.\hskip1pt B.\hskip1pt A.\hskip2pt Epstein, Curves on 2-manifolds and isotopies, {\it Acta Math.\ }115 (1966), 83--107.

[H1] A.\hskip1pt J.\hskip1pt S.\hskip2pt Hamilton, The triangulation of 3-manifolds, {\it Quart.\ J.\ Math.\ Oxford\/} 27 (1976), 63--70.

[H2] A.\hskip1pt Hatcher, {\it Algebraic Topology}, 2002.

[H3] A.\hskip1pt Hatcher, The Kirby torus trick for 3-manifolds, in preparation.

[K] R.\hskip1pt C.\hskip2pt Kirby, Stable homeomorphisms and the annulus conjecture, {\it Annals of Math.\ }89 (1969), 575--582.

[KS] R.\hskip1pt C.\hskip2pt Kirby and L.\hskip1pt C.\hskip2pt Siebenmann, {\it Foundational Essays on Topological Manifolds, Smoothings, and Triangulations}, Annals of Math. Studies 88, 1977.

[L] J.\hskip1pt M.\hskip1pt Lee, {\it Introduction to Smooth Manifolds}, Second Edition, 2013.

[M] E.\hskip1pt E.\hskip1pt Moise, {\it Geometric Topology in Dimensions 2 and 3}, 1977.

[R] T.\hskip1.5pt Rad\'o,   \"Uber den Begriff der Riemannschen Fl\"ache, {\it Acta Sci.\ Math.\ Szeged.\ }2 (1925), 101--121.

[T] C.\hskip1.5pt Thomassen,  The Jordan--Sch\"onflies theorem and the classification of surfaces, {\it Amer. Math. Monthly\ }99 (1992), 116--130.

\end